\documentclass[11pt,a4paper]{article}
\usepackage{amsfonts,amsmath,amssymb,accents,amsthm}%,mathtools}
\usepackage{verbatim}
\usepackage{hyperref}
\usepackage[normalem]{ulem}%for \sout,\xout

\usepackage{tikz} %also makes \color work
\newcommand{\dis}{\displaystyle}
\newcommand{\txt}{\textstyle}

\newcommand{\noi}{\noindent}
\newcommand{\halmos}{\rule{1ex}{1.4ex}}
\newcommand{\QED}{\nopagebreak{\hspace*{\fill}$\halmos$\medskip}}

\newcommand{\quand}{\quad\mbox{and}\quad}

\newtheoremstyle{mythm}% name
  {}%      Space above
  {}%      Space below
  {\itshape}%         Body font
  {}%         Indent amount (empty = no indent, \parindent = para indent)
  {\bfseries}% Thm head font
  {}%        Punctuation after thm head
  {.5em}%     Space after thm head: " " = normal interword space;
  {#1 #2 \thmnote{(#3)}}%  Thm head spec (can be left empty, meaning `normal')

\theoremstyle{mythm}
\newtheorem{theorem}{Theorem}[section]
\newtheorem{proposition}[theorem]{Proposition}
\newtheorem{lemma}[theorem]{Lemma}
\newtheorem{exercise}[theorem]{Exercise}
\newtheorem{corollary}[theorem]{Corollary}
\newtheorem{conjecture}[theorem]{Conjecture}

\newtheorem{counterex}[theorem]{Counterexample}
\newtheorem{remark}[theorem]{Remark}

\newcommand{\bt}{\begin{theorem}}
\newcommand{\et}{\end{theorem}}
\newcommand{\bl}{\begin{lemma}}
\newcommand{\el}{\end{lemma}}
\newcommand{\bp}{\begin{proposition}}
\newcommand{\ep}{\end{proposition}}
\newcommand{\bcor}{\begin{corollary}}
\newcommand{\ecor}{\end{corollary}}
\newcommand{\br}{\begin{remark}\rm}
\newcommand{\er}{\end{remark}}
\newcommand{\bcon}{\begin{conjecture}}
\newcommand{\econ}{\end{conjecture}}
\newcommand{\bex}{\begin{exercise}}
\newcommand{\eex}{\end{exercise}}
\newcommand{\bcou}{\begin{counterex}}
\newcommand{\ecou}{\end{counterex}}

\newenvironment{Proof}[1][]{\noi\textbf{Proof #1}}{\QED}
\newcommand{\bpro}{\begin{Proof}}
\newcommand{\epro}{\end{Proof}}

\newcommand{\be}{\begin{equation}}
\newcommand{\ee}{\end{equation}}
\newcommand{\ba}{\begin{array}}
\newcommand{\ea}{\end{array}}
\newcommand{\bc}{\be\begin{array}{r@{\,}c@{\,}l}}
\newcommand{\ec}{\end{array}\ee}

\newcommand{\de}{\delta}

\newcommand{\eps}{\varepsilon}
\newcommand{\la}{\lambda}

\newcommand{\sig}{\sigma}

\newcommand{\tet}{\theta}

\newcommand{\Ci}{{\cal C}}
\newcommand{\Di}{{\cal D}}

\newcommand{\Ki}{{\cal K}}

\newcommand{\Mi}{{\cal M}}
\newcommand{\Ni}{{\cal N}}

\newcommand{\R}{{\mathbb R}}
\newcommand{\N}{{\mathbb N}}
\newcommand{\Z}{{\mathbb Z}}

\newcommand{\E}{{\mathbb E}}
\renewcommand{\P}{{\mathbb P}}

\newcommand{\li}{\langle}
\newcommand{\re}{\rangle}

\newcommand{\volgt}{\ensuremath{\Rightarrow}}

\newcommand{\down}{\downarrow}
\newcommand{\sub}{\subset}

\newcommand{\toP}{\overset{\rm P}{\to}}

\newcommand{\Tofdd}{\overset{\rm f.d.d.}{\Longrightarrow}}

\newcommand{\asto}[1]{\underset{{#1}\to\infty}{\longrightarrow}}
\newcommand{\Asto}[1]{\underset{{#1}\to\infty}{\Longrightarrow}}

\newcommand{\ato}[2]{\underset{{#1}\to{#2}}{\longrightarrow}}
\newcommand{\Ato}[2]{\underset{{#1}\to{#2}}{\Longrightarrow}}
\newcommand{\atoP}[2]{\underset{{#1}\to{#2}}{\overset{\rm P}{\longrightarrow}}}
\newcommand{\Atofdd}[2]{\underset{{#1}\to{#2}}{\overset{\rm f.d.d.}{\Longrightarrow}}}

\newcommand{\ti}{\tilde}

\newcommand{\ov}{\overline}

\newcommand{\ffrac}[2]{{\textstyle\frac{{#1}}{{#2}}}}

\newcommand{\di}{\mathrm{d}}
\newcommand{\half}{{[0,\infty)}}

\newcommand{\ha}{\ffrac{1}{2}}

\newcommand{\var}{{\rm Var}}

\newcommand{\cadlag}{c\`adl\`ag}

\makeatletter
\let\@fnsymbol\@arabic
\makeatother

\begin{document}

%numbering formulas within sections
\makeatletter\@addtoreset{equation}{section}
\makeatother\def\theequation{\thesection.\arabic{equation}}

%alternative layout for enumerate lists.
\renewcommand{\labelenumi}{{\rm (\roman{enumi})}}
\renewcommand{\theenumi}{\roman{enumi}}

\title{An invariance principle for biased voter model interfaces}
\author{Rongfeng Sun \footnote{Department of Mathematics,
		National University of Singapore,
		10 Lower Kent Ridge Road, 119076 Singapore.
		Email: matsr@nus.edu.sg}
	\and Jan~M.~Swart \footnote{The Czech Academy of Sciences,
                Institute of Information Theory and Automation,
		Pod vod\'arenskou v\v e\v z\' i 4,
		18200 Praha 8,
		Czech Republic.
		Email: swart@utia.cas.cz}
	\and Jinjiong Yu \footnote{NYU-ECNU Institute of
		Mathematical Sciences at NYU Shanghai,
		3663 Zhongshan Road North,
		Shanghai 200062, China.
		Email: jinjiongyu@nyu.edu}
}

\date{July 24, 2020}

\maketitle

\begin{abstract}\noi
We consider one-dimensional biased voter models, where 1's replace 0's at a
faster rate than the other way round, started in a Heaviside initial state
describing the interface between two infinite populations of 0's and 1's. In
the limit of weak bias, for a diffusively rescaled process, we consider a
measure-valued process describing the local fraction of type 1 sites as a
function of time. Under a finite second moment condition on the rates, we show
that in the diffusive scaling limit there is a drifted Brownian path with the
property that all but a vanishingly small fraction of the sites on the left
(resp.\ right) of this path are of type 0 (resp.\ 1). This extends known
results for unbiased voter models. Our proofs depend crucially on recent
results about interface tightness for biased voter models.
\end{abstract}
\vspace{.5cm}

\noi
{\it Acknowledgements.} R.~Sun is supported by NUS grant R-146-000-253-114.
J.M.~Swart is supported by grant 19-07140S of the Czech Science Foundation (GA
CR).  J.~Yu is supported by a grant from NYU-ECNU Institute of Mathematical
Sciences at NYU Shanghai.

%82C22   	Interacting particle systems
%
%82C24 Interface problems; diffusion-limited aggregation
%82C41 Dynamics of random walks, random surfaces, lattice animals etc.
%60K35   	Interacting random processes; statistical mechanics type
%               models; percolation theory

{\setlength{\parskip}{-2pt}\tableofcontents}

\newpage

\section{Introduction}

\subsection{Statement of the result}

Let $\{0,1\}^\Z$ denote the space of all configurations of zeros and ones on
$\Z$, i.e., elements of $\{0,1\}^\Z$ are of the form $x=(x(i))_{i\in\Z}$ with
$x(i)\in\{0,1\}$. The one-dimensional biased voter model $(X^\eps_t)_{t\geq 0}$
with kernel $a(\,\cdot\,)$ and bias parameter $\eps\in[0,1)$ is the
  interacting particle system with state space $\{0,1\}^\Z$ and formal generator
\bc\label{bvmgen}
\dis G^\eps f(x)
&=&\dis\sum_{i,j}a(j-i)1_{\{x(i,j)=10\}}\big\{f(x+e_j)-f(x)\big\}\\[5pt]
&&\dis+(1-\eps)\sum_{i,j}a(j-i)1_{\{x(i,j)=01\}}\big\{f(x-e_j)-f(x)\big\},
\ec
where $e_i(j):=1_{\{i=j\}}$, and $x(i,j)=10$ is shorthand for $x(i)=1$,
$x(j)=0$. In words, (\ref{bvmgen}) says that if $x(i)=1$ and $x(j)=0$, then
the site $j$ adopts the type of site $i$ with rate $a(j-i)$. In the reverse case,
when $x(i)=0$ and $x(j)=1$, the site $j$ adopts the type of site $i$ with
rate $(1-\eps)a(j-i)$. In particular, for $\eps=0$, we obtain a standard voter model.

The kernel $a$ is a probability measure on $\Z$, not necessarily symmetric, such that $a(0)=0$.
In addition, throughout this paper, the following assumptions on $a$ will
always be in place:
\begin{enumerate}
\item $a$ is irreducible, i.e., each $k\in\Z$ can be written as a finite sum
  of $i\in\Z$ for which $a(i)>0$,
%\item $a$ is non-nearest-neighbor, i.e., $a(i)>0$ for some $i\in\Z$ with
%  $|i|>1$,
\item $a$ has mean zero, i.e., $\sum_ka(k)k=0$,
\item $a$ has a finite second moment, i.e., $\sig^2:=\sum_ka(k)k^2<\infty$.
%\be\label{sig}
%
%\ee
\end{enumerate}
We let
\be
S^{01}_{\rm int}:=\big\{x\in\{0,1\}^\Z:
\lim_{i\to-\infty}x(i)=0,\ \lim_{i\to\infty}x(i)=1\big\}
\ee
denote the space of states in which an infinite population of $0$'s on the
left and an infinite population of $1$'s on the right are separated by a
hybrid zone containing a mixture of $0$'s and $1$'s. This hybrid zone is called
the \emph{interface} of the biased voter model. If $X^\eps$ is started from an
initial state in $S^{01}_{\rm int}$ and $a$ has a finite first moment,
then it is known \cite{BMSV06} that almost surely $X^\eps_t\in S^{01}_{\rm
  int}$ for all $t\geq 0$.

If $a$ has a finite second moment, then it is moreover known \cite{SSY18} that
starting from a Heaviside state, the process returns to a Heaviside state
after a random time with finite expectation, so that the process modulo
translations is positive recurrent. This type of behavior is called
\emph{interface tightness} and was first proved for unbiased voter models in
the groundbreaking paper \cite{CD95} under a finite third moment condition on
the rates. The authors of \cite{BMV07} improved this to a finite second
  moment condition, which they showed is optimal.

Let $\Mi(\R)$ denote the space of locally finite measures on $\R$, equipped
with the topology of vague convergence. We use $(X^\eps_t)_{t\geq 0}$ to
define a measure-valued process $(\mu^\eps_t)_{t\geq 0}$ taking values in
$\Mi(\R)$ by
\be\label{mueps}
\mu^\eps_t:=\sum_{i\in\Z}\eps X^\eps_{\eps^{-2}t}(i)\de_{\eps i}
\qquad(t\geq 0),
\ee
where $\de_r$ denotes the delta-measure at $r\in\R$. We fix a standard
Brownian motion $(W_t)_{t\geq 0}$ and define a Brownian motion $B=(B_t)_{t\geq
  0}$ with drift $-\ha\sig^2$ and diffusion coefficient $\sig^2$ by
$B_t:=W_{\sig^2 t}-\ha\sig^2t$. We use $B$ to define a measure-valued
process $(\mu_t)_{t\geq 0}$ by
\be\label{mu}
\mu_t(\di x):=1_{\{x\geq B_t\}}\,\di x\qquad(t\geq 0),
\ee
i.e., $\mu_t$ has the density $1_{[B_t,\infty)}$ with respect to %w.r.t.\ 
the Lebesgue measure.
Our main result says that $(\mu_t)_{t\geq 0}$ arises as the weak limit of
$(\mu^\eps_t)_{t\geq 0}$.

\bt[Invariance principle for biased voter model interface]
Fix\label{T:measure} $x\in S^{01}_{\rm int}$ and for $\eps\in(0,1)$, let
$X^\eps$ be the biased voter model with generator \eqref{bvmgen} and initial
state $x$. Define $(\mu^\eps_t)_{t\geq 0}$ and $(\mu_t)_{t\geq 0}$ as in
\eqref{mueps} and \eqref{mu}. Then
\be\label{measure}
\P\big[(\mu^{\eps}_t)_{t\geq 0}\in\,\cdot\,\big]
\Ato{\eps}{0}\P\big[(\mu_t)_{t\geq 0}\in\,\cdot\,\big],
\ee
where $\Rightarrow$ denotes weak convergence
%of probability measures
on the Skorohod space $\Di\big([0,\infty),\Mi(\R)\big)$.
\et

\subsection{Main idea of the proof}\label{S:proid}

We will prove Theorem~\ref{T:measure} by first proving convergence of finite
dimensional distributions and then proving tightness. The convergence of
finite dimensional distributions will be obtained in a quick and elegant way
based on results we proved earlier in \cite{SSY18}. Our proof of tightness is
done via comparisons with the unbiased voter model, which still requires
significant work and new ideas.

For each $x\in S^{01}_{\rm int}$, there exists a unique
$M(x)\in\Z+\ha:=\{i+\ha:i\in\Z\}$ such that
\be
\sum_{i<M(x)}x(i)=\sum_{i>M(x)}\big(1-x(i)\big),
\ee
that is, the number of 1's to the left of the reference point $M(x)$ equals
the number of 0's to the right of it. We call $M(x)$ the \emph{weighted
  midpoint} of the interface. We call
\bc
\dis L(x)&:=&\dis\sup\big\{i\in\Z+\ha:x(j)=0\ \forall j<i\big\},\\[5pt]
\dis R(x)&:=&\dis\inf\big\{i\in\Z+\ha:x(j)=1\ \forall j>i\big\}
\ec
the \emph{left} and \emph{right boundary} of the interface, respectively.
Note that $L(x)=M(x)=R(x)$ if and only if $x$ is a \emph{Heaviside state} of the
form
\be
x_{{\rm hv},j}(i):=1_{\{i>j\}}\qquad(i\in\Z,\ j\in\Z+\ha).
\ee
In particular, we write $x_{\rm hv}:=x_{{\rm hv},1/2}$ for the state
  that has 1's on the positive integers and 0's elsewhere.
If $x$ is not a Heaviside state, then $L(x)<M(x)<R(x)$. 
As a first step towards proving
Theorem~\ref{T:measure}, we will prove the following weaker result.

\bt[Convergence of the weighted midpoint]
Fix\label{T:midconv} $x\in S^{01}_{\rm int}$ and for $\eps\in(0,1)$, let
$X^\eps$ be the biased voter model with generator \eqref{bvmgen} and initial
state $x$. Then
\be
\P\big[\big(\eps M(X^\eps_{\eps^{-2}t})\big)_{t\geq 0}\in\,\cdot\,\big]
\Ato{\eps}{0}\P\big[(B_t)_{t\geq 0}\in\,\cdot\,\big],
\ee
where $\Rightarrow$ denotes weak convergence
on the Skorohod space $\Di\big([0,\infty),\R\big)$, and $(B_t)_{t\geq 0}$ is
the drifted Brownian motion defined above \eqref{mu}.
\et

\noi
We remark that the idea of working on the weighted midpoint process is crucial.
On the one hand, in the scaling limit the left and right boundaries, 
and hence the measure-valued process, are well approximated by 
the weighted midpoint in finite dimensional distributions.
On the other hand, unlike the measure-valued process, there exists an efficient 
representation of the weighted midpoint, which is described below.

Whenever there is a site flipping from 1 (resp.\ 0) to 0 (resp.\ 1), the
weighted midpoint jumps to right (resp.\ left) by one.  Using this, it is easy
to see that the weighted midpoint evolves as a random time-changed random
walk, which has a drift of order $\eps$. In view of this, to prove
Theorem~\ref{T:midconv}, it suffices to control the random time
change. Somewhat surprisingly, it turns out that Lemma~3.1 and Proposition~3.7
of \cite{SSY18} give expressions for exactly the quantity we need and
Theorem~\ref{T:midconv} now follows from some relatively simple renewal
arguments. The paper \cite{SSY18} is concerned with \emph{interface
  tightness}, which we explain now.

We call two configurations $x,y\in\{0,1\}^\Z$ \emph{equivalent}, denoted by
$x\sim y$, if one is a translation of the other, i.e., there exists some
$k\in\Z$ such that $x(i)=y(i+k)$ $(i\in\Z)$. We let $\ov x$ denote the
equivalence class containing $x$ and write
\be\label{modtrans}
\ov S^{01}_{\rm int}:=\{\ov x:x\in S^{01}_{\rm int}\}.
\ee
Note that $S^{01}_{\rm int}$ and $\ov S^{01}_{\rm int}$ are countable sets. Since
our rates are translation invariant, the \emph{process modulo translations}
$(\ov{X}^\eps_t)_{t\geq0}$ is itself a Markov process; if we restrict the
state space to $\ov S^{01}_{\rm int}$, then it is in fact a continuous-time
Markov chain. If $a$ is non-nearest-neighbor, then it can be shown that this
Markov chain is irreducible (see \cite[Lemma~2.1]{SSY18}). 
In contrast, if $a$ is nearest-neighbor, then the Markov chain is absorbed at the 
Heaviside state $\ov x_{\rm hv}$ almost surely.
Following \cite{CD95}, we say that $(X^\eps_t)_{t\geq0}$ exhibits \emph{interface
  tightness on} $S^{01}_{\rm int}$ if $(\ov{X}^\eps_t)_{t\geq0}$ is positive 
 recurrent on $\ov S^{01}_{\rm int}$.
Under our assumptions (i)--(iii) on the kernel
$a$, interface tightness for biased voter models has been proved in
\cite[Thm~1.2]{SSY18}.

Interface tightness tells us that the biased voter model, started from any
initial state in $S^{01}_{\rm int}$, spends a positive fraction of its time in
Heaviside states. Moreover, the process modulo translations, started from
$\ov x_{\rm hv}$, returns to $\ov x_{\rm hv}$ in finite expected time.
Finally, the laws of the \emph{width} of the interface
$\P[R(X_t)-L(X_t)\in\,\cdot\,]$ are tight as $t\to\infty$. Theorem~1.3 of
\cite{SSY18} shows that all these statements hold uniformly as the bias $\eps$
tends to zero.

Combining \cite[Thms~1.2 and 1.3]{SSY18}, which prove interface tightness
uniformly as $\eps\down 0$, with the convergence of the weighted midpoint, we
then rather easily also obtain convergence in finite dimensional distributions
of the measure-valued process. To complete the proof of
Theorem~\ref{T:measure}, it therefore suffices to show tightness of the laws
of the measure-valued processes $(\mu^\eps_t)_{t\geq 0}$ as $\eps\down 0$. In
the unbiased setting, this has been proved in \cite{AS11} by directly verifying Jakubowski's tightness criterion (see e.g.\ \cite[Thm 3.6.4]{DA93}).
To use this criterion in the biased setting, we will construct a sufficient condition (\ref{tight1}) in Lemma~\ref{L:tight}, 
which requires to show that over a very short time interval, there cannot be too many 0's invading far into the region dominated by 1's and vice versa. 
The first scenario can be ruled out by a direct comparison with the unbiased voter model. 
The second scenario is more subtle and can be ruled out from the observation that, once enough 1's invade far into the region dominated by 0's, they will persist till the end of the short time interval, which will contradict the already proved fact that at any deterministic time, the interface location is distributed as that of a drifted Brownian motion. 
Combining the two scenarios, we are able to check (\ref{tight1}) and hence prove tightness.

\subsection{Discussion and open problems}\label{S:disc}

Combining Theorem~\ref{T:midconv} with \cite[Thm~1.3]{SSY18}, one can easily
show that as $\eps\downarrow0$, the diffusively rescaled left boundary
$\big(\eps L(X^\eps_{\eps^{-2}t})\big)_{t\geq0}$ and right boundary $\big(\eps
R(X^\eps_{\eps^{-2}t})\big)_{t\geq0}$ of the interface converge in finite
dimensional distributions to the same drifted Brownian motion as the weighted
midpoint. A natural question then arises. That is, as $\eps\downarrow0$, do
the boundaries also converge as processes, or equivalently does path level
tightness for the boundaries hold?

In the unbiased case $\eps=0$, this question has been answered in a sequence
of papers. Newman, Ravishankar and Sun \cite{NRS05} confirmed path level
tightness under the assumption that $a$ has a finite fifth moment. This result
was later extended by Belhaouari et al.\ in \cite{BMSV06} to all $a$ with a
finite $(3+\delta)$-th moment for some $\de>0$. On the other hand, it was pointed out in
\cite{BMSV06} that path level tightness for the left and right boundaries does
not hold if $\sum_{k} a(k)|k|^\gamma=\infty$ for some $\gamma<3$.

Indeed, in this regime, there exist exceptional times when $1$'s (resp.~$0$'s)
are created deep into the territory of the $1$'s (resp.~$0$'s) du
e to the
heavy tail of $a$. Nevertheless, such $1$'s and $0$'s are expected to be rare
and sparse, thanks to interface tightness. Therefore, one should be able to
restore tightness if those rare $1$'s and $0$'s are suitably discounted. In
\cite{BMSV06}, this idea was achieved by suppressing the infections $0\to1$
and $1\to0$ from site $i$ to site $j$ with $|i-j|\geq\epsilon^{-\kappa}$ for
some $\kappa>0$ depending on $a$, where $a$ is required to have a finite
$\gamma$-th moment for some $\gamma>2$.
The same idea also motivated Athreya and Sun \cite{AS11}, who proved an
unbiased version of Theorem~\ref{T:measure} assuming only that $a$ has a
finite second moment. 
We note that for interface tightness, the finite second moment assumption is optimal since 
it is shown in \cite{BMV07} that if $\sum_{k}a(k)|k|^\gamma=\infty$ for some $\gamma<2$, 
then interface tightness for the (unbiased) voter model does not hold. 

It is well-known that the voter model is dual to a system of coalescing random
walks. Likewise, the biased voter model is dual to a system of branching and
coalescing random walks. It has been shown in \cite{FINR04} that
nearest-neighbor systems of coalescing random walks, started from every point
in space and time, have a diffusive scaling limit, called the \emph{Brownian
  web}. Likewise, it has been shown in \cite{SS08} that nearest-neighbor
systems of weakly branching and coalescing random walks have a diffusive
scaling limit called the \emph{Brownian net}.
To extend these results to non-nearest-neighbor systems of (branching)
coalescing random walks, one needs to prove tightness for the collection of
paths in the Brownian web topology, introduced in \cite{FINR04}. In the
unbiased case, it has been shown in \cite{BMSV06} that tightness of coalescing
random walks in the Brownian web topology is equivalent to path level
tightness for the left and right boundaries of the dual voter model. Their
arguments carry over to branching coalescing random walks and their dual, the
biased voter model.

In view of this, for the biased voter model, it is an important open problem
to derive sufficient conditions for path level tightness for the left and
right boundaries. We conjecture that as in the unbiased case, a finite
$(3+\delta)$-th moment should suffice.

The remainder of the paper (which consists of Section~\ref{S:proofs} and an
appendix) is devoted to proofs. In Subsection~\ref{S:midconv} we give the main
line of the proof of Theorem~\ref{T:midconv} and in Subsections~\ref{S:time}
and \ref{S:renew} we fill in the details. In Subsections~\ref{S:fdd} and
\ref{S:tight} we then prove Theorem~\ref{T:measure} by first showing
convergence in finite dimensional distributions and then tightness. Lastly we collect 
some technical lemmas in the appendix.

\section{Proofs}\label{S:proofs}

\subsection{Convergence of the weighted midpoint}\label{S:midconv}

In this subsection we outline the proof of Theorem~\ref{T:midconv}. We show
that Theorem~\ref{T:midconv} follows from Lemmas~\ref{L:timechange} and
\ref{L:speed} below. Here Lemma~\ref{L:timechange} says that the weighted
midpoint evolves as a time-changed simple random walk, while
Lemma~\ref{L:speed} contains a statement about the convergence of the time
change. We also show how Lemma~\ref{L:speed} can heuristically be derived from
results proved in \cite{SSY18}, which prepares for its formal proof in
Subsection~\ref{S:renew} below.

For $x\in S^{01}_{\rm int}$ and $k\in\Z$, let
\be\label{Idef1}
I_k(x):=\{i:x(i)\neq x(i+k)\}
\ee
denote the number of \emph{$k$-boundaries} in the interface configuration $x$.
Note that $I_0(x)=0$ and $I_k(x)=I_{-k}(x)$ by definition.

\bl[Time-changed random walk]
Fix\label{L:timechange} $x\in S^{01}_{\rm int}$ and for $\eps\in(0,1)$, let
$X^\eps$ be the biased voter model with generator \eqref{bvmgen} and initial
state $x$. Then there exists an a.s.\ unique random, strictly increasing
continuous function $t\mapsto T^\eps_t$ such that
\be\label{timechange}
t=:\int_0^{T^\eps_t}\!\di s\,\sum_{k}a(k)I_k(X^\eps_s)\qquad(t\geq 0).
\ee
Moreover, the process $\big(M(X^\eps_{T^\eps_t})\big)_{t\geq 0}$ is a
continuous-time Markov chain on $\Z+\ha$ that jumps as
\be\label{driftrw}
m\mapsto m-1\mbox{ with rate }\ha
\quand
m\mapsto m+1\mbox{ with rate }\ha(1-\eps).
\ee
\el

By standard results, the drifted random walk in (\ref{driftrw}) converges
after diffusive rescaling to a drifted Brownian motion. In view of this, in
order to prove Theorem~\ref{T:midconv}, the main task is to control the
time-change in (\ref{timechange}). 
Interestingly, the time-change expression also appeared in the 
voter model equilibrium equation (see (\ref{equil}) below), from which we
can show that the time-change converges uniformly
to a deterministic limit. 

\bl[Convergence of the time change]
Let\label{L:speed} $X^\eps$ be as in Theorem~\ref{T:midconv}. Then
\be\label{lln1}
\sup_{0\leq t\leq T}
\Big|\sig^2t-\eps^{2}\int_0^{\eps^{-2}t}\!\di s\,\sum_ka(k)I_k(X^\eps_s)\Big|
\atoP{\eps}{0}0\qquad(T<\infty),
\ee
where $\toP$ denotes convergence in probability.
\el

\bpro[Proof~of Theorem~\ref{T:midconv}]
Let us set
\be
Y^\eps_t:=\eps M(X^\eps_{T^\eps_{\eps^{-2}t}})\qquad(t\geq 0),
\ee
i.e., this is the drifted random walk $(M(X^\eps_{T^\eps_t}))_{t\geq 0}$ from
Lemma~\ref{L:timechange}, diffusively rescaled by $\eps$. Then standard results
tell us that
\be
\P\big[(Y^\eps_{t})_{t\geq 0}\in\,\cdot\,\big]
\Ato{\eps}{0}\P\big[(W_t-\ha t)_{t\geq 0}\in\,\cdot\,\big],
\ee
where $(W_t)_{t\geq 0}$ is a standard Brownian motion. Let
\be\label{Tmin}
S^\eps_t=\int_0^t\!\di s\,\sum_ka(k)I_k(X^\eps_s)\qquad(t\geq 0)
\ee
denote the inverse of the function $t\mapsto T^\eps_t$ defined in
(\ref{timechange}). Then
\be\label{YT}
\eps M(X^\eps_{\eps^{-2}t})=Y^\eps_{\eps^2S^\eps_{\eps^{-2}t}}\qquad(t\geq 0).
\ee
Lemma~\ref{L:speed} tells us that $\eps^2 S^\eps_{\eps^{-2}t}$ converges as a
process to $\sig^2 t$. It is not hard to show (for details we refer to
Lemma~\ref{L:timconv} in the appendix) that this implies convergence of the
time-changed process, proving the claim of Theorem~\ref{T:midconv}.
\epro

In order to prove the crucial Lemma~\ref{L:speed}, we will rely on results
proved in \cite{SSY18}. In the remainder of this subsection, we recall some of
these results and put them into context, to give the reader a rough idea where
Lemma~\ref{L:speed} comes from.

As explained in Subsection~\ref{S:proid}, Theorem~1.2 in \cite{SSY18}
establishes interface tightness for biased voter models. More precisely, this
theorem says that for any $\eps\in[0,1)$, the process modulo translations
$(\ov{X}^\eps_t)_{t\geq0}$ has a unique invariant law. Let us denote this
invariant law by $\ov\pi^\eps$. In particular, if $a$ is nearest-neighbor (and
therefore $a(-1)=\ha=a(1)$ by our assumptions on $a$), then the Heaviside
state $\ov x_{\rm hv}$ is absorbing and $\ov\pi^\eps$ is the delta measure on
$\ov x_{\rm hv}$.  In the non-nearest-neighbor case, we cite the following
theorem from \cite[Thm~1.3]{SSY18}. The extension to the nearest-neighbor case
is trivial.

\bt[Continuity of the invariant law]
The\label{T:invlaw} laws $\ov\pi^\eps$ converge weakly to $\ov\pi^0$ as
$\eps\downarrow0$ with respect to the discrete topology on $\ov S^{01}_{\rm
  int}$.
\et

All existing proofs of interface tightness for unbiased voter models are in
some way or another based on a function that counts the number of
\emph{inversions}, i.e., pairs of sites $i,j$ such that $i<j$ and
$x(i)>x(j)$. Let $h$ denote this function, i.e.,
\be\label{numinv}
h(x):=\sum_{i<j}1_{\txt\{x(i,j)=10\}}\qquad(x\in S^{01}_{\rm int}).
\ee
Note that since $h$ is translation invariant, we can alternatively view $h$ as
a function on $\ov S^{01}_{\rm int}$. In \cite[Prop.~3.7]{SSY18}, it is shown
that the invariant law $\ov\pi^0$ solves the equilibrium equation
\be
\sum_{\ov x\in\ov S^{01}_{\rm int}}\ov\pi^0(\ov x)\,G^0h(\ov x)=0,
\ee
where $G^0$ is the generator defined in (\ref{bvmgen}). As shown in
\cite[Prop.~3.7]{SSY18}, this equation can be written more explicitly as
follows. (In the nearest-neighbor case, \cite[Prop.~3.7]{SSY18} is not
applicable, but (\ref{equil}) below holds trivially with both sides equal to
1.)

\bp[Equilibrium equation]
Let\label{P:equil} $X^0_\infty$ be a random variable such that $\ov X^0_\infty$
has law $\ov\pi^0$. Then
\be\label{equil}
\E\Big[\sum_ka(k)I_k(X^0_\infty)\Big]=\sig^2.
\ee
\ep

It is a remarkable fact that the equilibrium equation for the function in
(\ref{numinv}) yields an expression for precisely the quantity that also
appears in the time-change in (\ref{timechange}). Proposition~\ref{P:equil}
was one of the main ingredients used in \cite{SSY18} to prove
Theorem~\ref{T:invlaw}.  Together, Theorem~\ref{T:invlaw} and
Proposition~\ref{P:equil} will be the main ingredients in our proof of
Lemma~\ref{L:speed}. In order to derive Lemma~\ref{L:speed} from
(\ref{equil}), we will need uniform control over the speed at which the
process modulo translations converges to equilibrium. This will be achieved by
a renewal decomposition of the process modulo translations, where we get
uniform control on the expected return times as $\eps\down 0$ as a result of
Theorem~\ref{T:invlaw}.

In the coming two subsections, we prove Lemmas~\ref{L:timechange} and
\ref{L:speed}, respectively.

\subsection{Random time-change}\label{S:time}

\bpro[Proof~of Lemma~\ref{L:timechange}]
Let
\be\label{Idef2}
\ba{l}
I^{01}_k(x):=\{i:x(i)=0,x(i+k)=1\},\\ [5pt]
I^{10}_k(x):=\{i:x(i)=1,x(i+k)=0\}.
\ec
Since for every $x\in S^{01}_{\rm int}$, $k>0$ and $i\in\Z$, there is one more adjacent pair of (01) than (10) along the subsequence $x(\ldots,i-k,i,i+k,\ldots)$, it is not hard to see that
\be
I^{01}_k(x)=I^{10}_k(x)+k\quand
I_k(x)=I^{01}_k(x)+I^{10}_k(x)\qquad(x\in S^{01}_{\rm int}).
\ee
As a result
\be\label{Iref2}
I^{10}_k(x)=\ha(I_k(x)-k)\quand I^{01}_k(x)=\ha(I_k(x)+k).
\ee
We observe that the quantity $M(X^\eps_t)$ always goes up and down by a single
unit. More precisely, $M(X^\eps_t)$ goes down by one when a site flips from
$0$ to $1$ and it goes up by one when a site flips from $1$ to $0$, which
means that if the present state is $X^\eps_t=x$, then $M(X^\eps_t)$ jumps as
\be\ba{r@{\quad}c@{\quad}l}\label{Midrates}
m\mapsto m-1&\mbox{with rate}
&\dis\sum_ka(k)I^{10}_k(x)=\ha\sum_ka(k)I_k(x),\\[5pt]
m\mapsto m+1&\mbox{with rate}
&\dis(1-\eps)\sum_ka(k)I^{01}_k(x)=(1-\eps)\ha\sum_ka(k)I_k(x),
\ec
where we have used (\ref{Iref2}) and our assumption $\sum_ka(k)k=0$.
It follows from (\ref{Midrates}) that $M(X^\eps_t)$ is a drifted random walk with a random time change.

More precisely, defining $S^\eps_t$ as in (\ref{Tmin}), and observing that the
integrand is $\geq 1$, we see that $S^\eps_t$ is a.s.\ strictly increasing,
continuous, with $S^\eps_0=0$ and $\lim_{t\to\infty}S^\eps_t=\infty$. It
follows that $S^\eps_t$ has an inverse function with the same properties,
which is $T^\eps_t$. By standard results, the time-changed process
$(X^\eps_{T^\eps_t})_{t\geq 0}$ is a Markov process such that if the original
process $(X^\eps_t)_{t\geq 0}$ jumps from $x$ to $y$ with rate $r(x,y)$, then
the new process $(X^\eps_{T^\eps_t})_{t\geq 0}$ jumps from $x$ to $y$ with rate
$(\sum_ka(k)I_k(x))^{-1}r(x,y)$. In particular, the process
$\big(M(X^\eps_{T^\eps_t})\big)_{t\geq 0}$ is a drifted random walk with jump rates
as in (\ref{driftrw}).
\epro

\subsection{Renewal arguments}\label{S:renew}

In this subsection, we prove Lemma~\ref{L:speed}, completing the proof of
Theorem~\ref{T:midconv}. Since the functions $I_k$ are translation invariant,
we can and will view them as functions on $\ov S^{01}_{\rm int}$. Our task is
then to show that if $\ov x\in \ov S^{01}_{\rm int}$ is fixed and
$\ov X^\eps$ is the biased voter model with bias $\eps\in(0,1)$, modulo
translations, started in $\ov x$, then
\be\label{ovlln1}
\sup_{0\leq t\leq T}
\Big|\sig^2t-\eps^{2}\int_0^{\eps^{-2}t}\!\di s\,\sum_ka(k)I_k(\ov X^\eps_s)\Big|
\atoP{\eps}{0}0.
\ee
We let
\be\label{tau0}
\tau^\eps_0:=\inf\{t\geq 0:\ov X^\eps_t=\ov x_{\rm hv}\}
\ee
denote the first hitting time of $\ov x_{\rm hv}$, and define inductively
\be\label{retime}
\ti\tau^\eps_n:=\inf\big\{t>\tau^\eps_{n-1}:\ov X^\eps_t\neq\ov x_{\rm hv}\big\}
\quand
\tau^\eps_n:=\inf\big\{t>\ti\tau^\eps_n:\ov X^\eps_t=\ov x_{\rm hv}\big\}
\qquad(n\geq 1).
\ee
We will first prove (\ref{ovlln1}) under the additional assumptions that
$\ov X^\eps_0=\ov x_{\rm hv}$ and the kernel $a$ is non-nearest-neighbor.
The assumption that $\ov X^\eps_0=\ov x_{\rm hv}$ implies that $\tau^\eps_0$
from (\ref{tau0}) is zero, while the assumption that $a$ is non-nearest-neighbor
implies that $r_\eps>0$, where
\be\label{reps}
r_\eps:=\sum_{k<-1}(|k|-1)a(k)+(1-\eps)\sum_{k>1}(k-1)a(k)
\ee
is the rate at which $\ov X^\eps$ jumps away from $\ov x_{\rm hv}$. We start
with a trivial observation. Below, we view the law of $\ov X^\eps$ as a
probability measure on the space of piecewise constant, right-continuous
functions with values in the countable set $\ov S^{01}_{\rm int}$, and we
equip this space with the Skorohod topology.

\bl[Continuity of the law]
Let\label{L:colaw} $\ov X^\eps=(\ov X^\eps_t)_{t\geq 0}$ be the biased voter
model modulo translations with bias $\eps\in[0,1)$, started in
$\ov X^\eps_0=\ov x_{\rm hv}$. Then the function $\eps\mapsto\P\big[\ov
    X^\eps\in\,\cdot\,]$ is continuous with respect to weak convergence.
\el

\bpro
This is trivial, since $\ov X^\eps$ is a nonexplosive continuous-time Markov
chain and its jump rates converge pointwise.
\epro

\bl[Convergence of return times]
Assume\label{L:retconv} that $a$ is non-nearest-neighbor. Then
\be\label{retbdd}
\lim_{\eps\downarrow0}\E^{\ov x_{\rm hv}}[\tau^\eps_1]
=\E^{\ov x_{\rm hv}}[\tau^0_1]<\infty.
\ee
\el

\bpro
By interface tightness
\cite[Thm~1.2]{SSY18}, we have $\E[\tau^\eps_1]<\infty$ for each $\eps\in[0,1)$.
The regenerative theorem (see \cite[Thm~4.1.2]{A03}) gives
an expression for the invariant law $\ov\pi^\eps$,
\be\label{invlaw}
\ov\pi^\eps(\ov{x})=\frac{1}{\E^{\ov x_{\rm hv}}[\tau^\eps_1]}
\E^{\ov x_{\rm hv}}\Big[\int_{0}^{\tau^\eps_1}1_{\{\ov{X}^\eps_s=\ov{x} \}}\di s\Big]
\qquad(\ov{x}\in\ov{S}^{01}_{\rm int}).
\ee
In particular, setting $\ov{x}=\ov{x}_{\rm hv}$ and using the fact that
  during $[0,\tau^\eps_1)$ one has $\ov X^\eps_s=\ov x_{\rm hv}$ if and only if $s\in[0,\ti\tau^\eps_1)$, it follows that
\be\label{invpi}
\E^{\ov x_{\rm hv}}[\tau^\eps_1]
=\frac{1}{\ov\pi^\eps(\ov x_{\rm hv})}\E^{\ov x_{\rm hv}}[\ti\tau^\eps_1]
=\frac{1}{r_\eps\ov\pi^\eps(\ov x_{\rm hv})},
\ee
where $r_\eps$ from (\ref{reps}) is the rate at which $\ov X^\eps$ leaves $\ov
x_{\rm hv}$. By Theorem~\ref{T:invlaw}, $\ov\pi^\eps(\ov x_{\rm
  hv})\to\ov\pi^0(\ov x_{\rm hv})$ as $\eps\to 0$, which together with
(\ref{invpi}) yields the claim.
\epro

\bl[Average value during one excursion]
Assume\label{L:avconv} that $a$ is non-nearest-neighbor. Then
\be\label{avbdd}
\lim_{\eps\downarrow0}
\E^{\ov x_{\rm hv}}\Big[\int_0^{\tau^\eps_1}\!\di s\,\sum_ka(k)I_k(\ov X^\eps_s)\Big]
=\E^{\ov x_{\rm hv}}\Big[\int_0^{\tau^0_1}\!\di s\,\sum_ka(k)I_k(\ov X^0_s)\Big]
<\infty.
\ee
\el

\bpro
Formula (\ref{invlaw}) gives
\be\label{excexpr}
\E^{\ov x_{\rm hv}}\Big[\int_0^{\tau^\eps_1}\!\di s\,\sum_ka(k)I_k(\ov X^\eps_s)\Big]
=\E^{\ov x_{\rm hv}}[\tau^\eps_1]
\sum_{\ov x\in\ov S^{01}_{\rm int}}\ov\pi^\eps(\ov x)\sum_ka(k)I_k(\ov x).
\ee
By Theorem~\ref{T:invlaw} and Fatou's lemma
\be\label{piFat}
\sum_{\ov x\in\ov S^{01}_{\rm int}}\ov\pi^0(\ov{x})\sum_{k}a(k)I_k(\ov x)
\leq\liminf_{\eps\down 0}\sum_{\ov x\in\ov S^{01}_{\rm int}}
\ov\pi^\eps(\ov{x})\sum_ka(k)I_k(\ov x).
\ee
By \cite[Lemma~3.1]{SSY18} and Proposition~\ref{P:equil},
\be\label{expV2}
\sum_{\ov x\in\ov S^{01}_{\rm int}}\ov\pi^\eps(\ov{x})\sum_{k}a(k)I_k(\ov x)\leq\sig^2
\quand
\sum_{\ov x\in\ov S^{01}_{\rm int}}\ov\pi^0(\ov{x})\sum_{k}a(k)I_k(\ov x)=\sig^2.
\ee
The right-hand side of (\ref{piFat}) is less than or equal to the limit superior
of the same expression, which by the first formula in (\ref{expV2}) is 
bounded from above by $\sig^2$, while the second formula in (\ref{expV2})
identifies the left-hand side of (\ref{piFat}) as $\sig^2$. We conclude that
\be
\sum_{\ov x\in\ov S^{01}_{\rm int}}\ov\pi^\eps(\ov{x})\sum_ka(k)I_k(\ov x)
\ato{\eps}{0}\sig^2
=\sum_{\ov x\in\ov S^{01}_{\rm int}}\ov\pi^0(\ov{x})\sum_{k}a(k)I_k(\ov x).
\ee
Inserting this into (\ref{excexpr}), using Lemma~\ref{L:retconv}, we obtain
(\ref{avbdd}).
\epro

The proof of Lemma~\ref{L:avconv} yields a useful corollary.

\bcor[Renewal identity]
Assume\label{C:regen} that $a$ is non-nearest-neighbor. Then
\be
\frac{1}{\E^{\ov x_{\rm hv}}[\tau^0_1]}\,
\E^{\ov x_{\rm hv}}\Big[\int_0^{\tau^0_1}\!\di s\,\sum_ka(k)I_k(\ov X^0_s)\Big]
=\sig^2.
\ee
\ecor

\bpro
This follows from (\ref{excexpr}) and (\ref{expV2}).
\epro

Let $\ov X^\eps$ denote the process started in $\ov X_0^\eps=\ov x_{\rm hv}$ and
let
\be\label{phipsidef}
\phi_\eps(u):=\eps^2\sum_{k=1}^{\lfloor\eps^{-2}u\rfloor}(\tau^\eps_k-\tau^\eps_{k-1})
\quand
\psi_\eps(u):=\eps^2\sum_{k=1}^{\lfloor\eps^{-2}u\rfloor}
\int_{\tau^\eps_{k-1}}^{\tau^\eps_k}\!\di s\,\sum_ka(k)I_k(\ov X^\eps_s).
\ee
By the strong Markov property, $\phi_\eps(u)$ and $\psi_\eps(u)$ are sums of i.i.d.\ random
variables. Indeed, $\tau^\eps_k-\tau^\eps_{k-1}$ is equally distributed with
$\tau^\eps_1$ while the summands of $\psi_\eps(u)$ are equally distributed
with
\be
\eta^\eps:=\int_0^{\tau^\eps_1}\!\di s\,\sum_ka(k)I_k(\ov X^\eps_s).
\ee
It follows from Lemma~\ref{L:colaw} that $\tau^\eps_1$ and $\eta^\eps$
converge weakly in law as $\eps\to 0$ to $\tau^0_1$ and $\eta^0$,
respectively. Note that $\tau^0_1,\eta^0>0$ a.s. Lemmas~\ref{L:retconv},
\ref{L:avconv}, and Corollary~\ref{C:regen} tell us that
\be\label{meanconv}
\lim_{\eps\to 0}\E[\tau^\eps_1]=\E[\tau^0_1]<\infty,\quad
\lim_{\eps\to 0}\E[\eta^\eps]=\E[\eta^0]<\infty,\quand
\E[\eta^0]/\E[\tau^0_1]=\sig^2.
\ee
By \cite[Prop.~A.2.3]{EK86}, the convergence in law and in expectation of
$\tau^\eps_1$ and $\eta^\eps$ imply that these random variables are uniformly
integrable as $\eps\down 0$, i.e., for any $\eps_n\to 0$, we have
\be
\lim_{K\to\infty}\sup_n\E[\tau^{\eps_n}_1;\tau^{\eps_n}_1\geq K]=0
\quand
\lim_{K\to\infty}\sup_n\E[\eta^{\eps_n};\eta^{\eps_n}\geq K]=0.
\ee
It follows that
\be
\lim_{\eps\down 0}\E[\tau^\eps_1;\tau^\eps_1>t\eps^{-2}]=0
\quand
\lim_{\eps\down 0}\E[\eta^\eps;\eta^\eps>t\eps^{-2}]=0\qquad(t>0).
\ee
This allows us to apply a standard functional law of large numbers (see
Lemma~\ref{L:FLLN} in the appendix for details) to obtain the following lemma.

\bl[Functional law of large numbers]
Let\label{L:phipsi} $\phi_\eps$ and $\psi_\eps$ be as in
(\ref{phipsidef}). Then
\be
\sup_{0\leq u\leq U}\big|u\E[\tau^0_1]-\phi_\eps(u)\big|\atoP{\eps}{0}0
\quand
\sup_{0\leq u\leq U}\big|u\E[\eta^0]-\psi_\eps(u)\big|\atoP{\eps}{0}0\qquad(U>0).
\ee
\el

\bpro[Proof~of Lemma~\ref{L:speed}]
We will finish the proof in two steps, first assuming a Heaviside initial state 
and then extending to arbitrary initial states.\\
{\bf Heaviside initial state.}
We first prove the statement under the additional assumptions that
the kernel $a$ is non-nearest-neighbor and $\ov X^\eps_0=\ov x_{\rm hv}$.
In this case $\tau^\eps_0=0$.

Since $\phi_\eps(\eps^2 k)=\eps^2\tau^\eps_k$ $(k\geq 0)$, the function
$\phi_\eps$ defines a bijection from $\eps^2\N$ to $\{\eps^2\tau^\eps_k:k\geq
0\}$. Let $\tet_\eps$ denote the restriction of $\phi_\eps$ to
$\eps^2\N$. Then $\phi_\eps$ is the right-continuous interpolation of
$\tet_\eps$. Let $\tet^{-1}_\eps$ denote the inverse of $\tet_\eps$.
For any $t\geq 0$ and $k\in\N$, define
\be
[t]^\eps_\leftarrow:=\tau^\eps_{k-1}\quad\mbox{where}\quad
t\in[\tau^\eps_{k-1},\tau^\eps_k)\quand
[t]^\eps_\rightarrow:=\tau^\eps_k\quad\mbox{where}\quad
t\in(\tau^\eps_{k-1},\tau^\eps_k].
\ee
Then for $t\geq0$,
\be\label{twobd}
\psi_\eps\big(\tet^{-1}_\eps(\eps^2[\eps^{-2}t]^\eps_\leftarrow)\big)
\leq\eps^{2}\int_0^{\eps^{-2}t}\!\di s\,\sum_ka(k)I_k(\ov X^\eps_s)
\leq\psi_\eps\big(\tet^{-1}_\eps(\eps^2[\eps^{-2}t]^\eps_\rightarrow)\big),
%\qquad(t\geq 0),
\ee
where $\tet^{-1}_\eps(\eps^2[\eps^{-2}t]^\eps_\leftarrow)$ and
$\tet^{-1}_\eps(\eps^2[\eps^{-2}t]^\eps_\rightarrow)$ are the right- and
left-continuous interpolations of the function $\tet^{-1}_\eps$.

Lemma~\ref{L:phipsi} tells us
that as $\eps\down 0$, the right-continuous interpolation of $\tet_\eps$
converges in probability w.r.t.\ the Skorohod topology to the function
$u\mapsto u\E[\tau^0_1]$. By the Skorohod representation theorem
\cite[Thm~6.7]{Bil99}, along any sequence $\eps_n\down 0$, we can couple our
random variables such that this convergence is a.s. Since $\tet_\eps$
takes values in $\{\eps^2\tau^\eps_k:k\geq 0\}$, it is easy to see that
for our coupling
\be\label{denstim}
\forall t\geq 0\quad\exists t_n\in\{\eps_n^2\tau^{\eps_n}_k:k\geq 0\}
\quad\mbox{s.t.}\quad t_n\to t,
\ee
i.e., the range of $\tet_{\eps_n}$ is a.s.\ dense in the limit. Since the
sets $\eps_n^2\N$ are dense in the limit, it is easy to see that not only the
right-continuous interpolation, but also the linear interpolation of
$\tet_{\eps_n}$ converges locally uniformly to the function $u\mapsto
u\E[\tau^0_1]$. It is not hard to see that this implies locally uniform
convergence of the inverse (see Lemma~\ref{L:invconv} in the appendix).  Thus,
the linear interpolation of $\tet^{-1}_{\eps_n}$ converges locally uniformly
to the function $t\mapsto t/\E[\tau^0_1]$. Using (\ref{denstim}), we see that
the same holds for the right- and left-continuous interpolations of the
function $\tet^{-1}_{\eps_n}$. Since this holds for arbitrary $\eps_n\down 0$,
we obtain that
\be
\sup_{0\leq t\leq T}\big|t/\E[\tau^0_1]
-\tet^{-1}_\eps(\eps^2[\eps^{-2}t]^\eps_\leftarrow)\big|\atoP{\eps}{0}0,
\ee
and similarly for the left-continuous interpolation. Combining this with
Lemma~\ref{L:phipsi}, it is not hard to show (see Lemma~\ref{L:timconv} from
the appendix) that the left- and right-hand sides of (\ref{twobd}) converge
locally uniformly in probability to the composition of the functions $t\mapsto
t/\E[\tau^0_1]$ and $u\mapsto u\E[\eta^0]$. By (\ref{meanconv}), this
composite function is the function $t\mapsto\sig^2t$, proving (\ref{lln1}).

This completes the proof under the additional assumptions that
the kernel $a$ is non-nearest-neigbor and $\ov X^\eps_0=\ov x_{\rm hv}$.
If $a$ is the nearest-neighbor kernel $a(-1)=\ha=a(1)$, then $\sig^2=1$ and
$\ov X^\eps_t=\ov x_{\rm hv}$ for each $t\geq 0$. Moreover,
$\sum_ka(k)I_k(\ov x_{\rm hv})=1$, so in this case (\ref{ovlln1}) is trivial.

\noindent {\bf Arbitrary initial states.} To treat the case when $\ov X^\eps$ started in an arbitrary, fixed initial
state $\ov X^\eps_0=\ov x$, it suffices to show that
\be\label{hvtrap}
\eps^2\tau^\eps_0\atoP{\eps}{0}0\quand
\eps^2\int_0^{\tau^\eps_0}\!\di s\,\sum_ka(k)I_k(\ov X^\eps_s)\atoP{\eps}{0}0.
\ee
Since the jump rates converge, for any $\ov x\in\ov S^{01}_{\rm int}$, the laws
\be\label{asfin}
\P^{\ov x}\big[\tau^\eps_0\in\,\cdot\,\big]
\quand
\P^{\ov x}\big[\int_0^{\tau^\eps_0}\!\di s\,
\sum_ka(k)I_k(\ov X^\eps_s)\in\,\cdot\,\big]
\ee
converge weakly as $\eps\down 0$, so it suffices to show that for the unbiased
process
\be
\P^{\ov x}\big[\tau^0_0<\infty\big]=1
\quand
\P^{\ov x}\big[\int_0^{\tau^0_0}\!\di s\,\sum_ka(k)I_k(\ov X^0_s)<\infty\big]=1,
\ee
where in fact the first equality implies the second equality. If the kernel
$a$ is non-nearest-neighbor, then the fact that $\tau^0_0<\infty$ for any $\ov
x\in\ov S^{01}_{\rm int}$ follows from positive recurrence and
irreducibility \cite[Thm~1.2 and Lemma~2.1]{SSY18}.

To complete the proof, we must show that the nearest-neighbor unbiased voter
model, started in any state $x\in S^{01}_{\rm int}$, a.s.\ reaches a Heaviside
state in finite time. We can obtain a two-type voter model as a function of a
multi-type voter model in which initially each site has a different
type. Since the family size of each type is a martingale, each family dies out
a.s. As soon as the families corresponding to types that were initially
between $L(x)$ and $R(x)$ have all died out, $X_t$ and thus $\ov X_t$ will be in a Heaviside
state.
\epro

\subsection{Convergence of finite dimensional distributions}\label{S:fdd}

In this subsection, we start proving Theorem~\ref{T:measure} by showing
convergence in finite dimensional distributions.

\bl[Local limit]
Fix\label{L:loclim} $\ov x\in\ov S^{01}_{\rm int}$ and for $\eps\in[0,1)$, let
$\ov X^\eps$ be the biased voter model modulo translations with bias $\eps$,
started in $\ov x$. Then, for each $\eps_n\to 0$ and $t_n\to\infty$,
\be\label{loclim}
\P^{\ov x}\big[\ov X^{\eps_n}_{t_n}\in\,\cdot\,\big]\Asto{n}\ov\pi^0,
\ee
where $\Rightarrow$ denotes weak convergence of probability measures on $\ov
S^{01}_{\rm int}$ with respect to the discrete topology, and $\ov\pi^0$ is the
invariant law of $\ov X^0$.
\el

\bpro
We first prove the statement if the kernel $a$ is non-nearest-neighbor. The
process $\ov X^\eps$ is irredicible (by \cite[Lemma~2.1]{SSY18}) and positive recurrent (by \cite[Thm~1.2]{SSY18}) for each $\eps\geq
0$. Moreover, its jump rates and
by Theorem~\ref{T:invlaw} also its invariant law converge to those of $\ov X^0$ as
$\eps\down 0$. Using this, a simple abstract argument (see Lemma~\ref{L:unerg}
in the appendix) gives
\be
\sup_{n\geq 0}
\big\|\P^{\ov x}[\ov X^{\eps_n}_t\in\,\cdot\,]-\ov\pi^{\eps_n}\big\|
\asto{t}0,
\ee
where $\|\,\cdot\,\|$ denotes the total variation norm. Since
\be
\big\|\P^{\ov x}\big[\ov X^{\eps_n}_{t_n}\in\,\cdot\,\big]-\ov\pi^0\big\|
\leq
\big\|\P^{\ov x}\big[\ov X^{\eps_n}_{t_n}\in\,\cdot\,\big]-\ov\pi^{\eps_n}\big\|
+\big\|\ov\pi^{\eps_n}-\ov\pi^0\big\|,
\ee
the claim follows from the convergence of $\ov\pi^{\eps_n}$
(Theorem~\ref{T:invlaw}).

Since for the nearest-neighbor kernel $a(-1)=\ha=a(1)$, the invariant
law $\ov\pi^0$ is the delta measure on $\ov x_{\rm hv}$, in this case
it suffices to prove that
\be\label{diag0}
\P^{\ov x}\big[\tau^{\eps_n}_0>t_n]\asto{n}0,
\ee
where $\tau^{\eps_n}_0$ as in (\ref{tau0}) denotes the time $\ov X^{\eps_n}$ gets
trapped in $\ov x_{\rm hv}$. It has already been shown below (\ref{asfin}) that
\be
\lim_{\eps\down 0}\P^{\ov x}\big[\tau^\eps_0>t]=\P^{\ov x}\big[\tau^0_0>t]
\quand
\lim_{t\to\infty}\P^{\ov x}\big[\tau^0_0>t]=0,
\ee
which implies (\ref{diag0}).
\epro

\bp[Convergence of the left and right boundaries]
Fix\label{P:LRfdd} $x\in S^{01}_{\rm int}$ and for $\eps\in(0,1)$, let
$X^\eps$ be the biased voter model with generator (\ref{bvmgen}) and initial
state $x$. Then
\be\label{LRfdd}
\P\big[
\big(\eps L(X^\eps_{\eps^{-2}t}),\eps R(X^\eps_{\eps^{-2}t})\big)_{t\geq 0}
\in\,\cdot\,\big]
\Atofdd{\eps}{0}\P\big[(B_t,B_t)_{t\geq 0}\in\,\cdot\,\big],
\ee
where $\Tofdd$ denotes weak convergence of the finite dimensional
distributions, and $(B_t)_{t\geq 0}$ is a Brownian motion with drift
$-\ha\sig^2$ and diffusion coefficient $\sig^2$.
\ep

\bpro
Let $W(x):=R(x)-L(x)$ denote the width of an interface $x\in S^{01}_{\rm
  int}$, which by translation invariance we can view as a function on $\ov
S^{01}_{\rm int}$. Lemma~\ref{L:loclim} shows that as $\eps\down 0$, the law of
$W(\ov X^\eps_{\eps^{-2}t})$ converges to a limit law on $\N$, and hence
\be
\P^x\big[\eps W(X^\eps_{\eps^{-2}t})\in\,\cdot\,\big]\Ato{\eps}{0}0
\qquad(t>0).
\ee
Since $x$ is fixed, this trivially also holds for $t=0$. By
Theorem~\ref{T:midconv} and the Skorohod representation theorem
\cite[Thm~6.7]{Bil99}, along any sequence $\eps_n\down 0$,
we can couple our processes such that $(\eps_n
M(X^{\eps_n}_{\eps_nn^{-2}t}))_{t\geq 0}$ converges a.s.\ to $(B_t)_{t\geq 0}$.
The claim now follows since $|L(x)-M(x)|\leq W(x)$ and similarly for $R(x)$
$(x\in S^{01}_{\rm int})$.
\epro

\bl[Convergence of finite dimensional distributions]
Fix\label{L:fdd} $x\in S^{01}_{\rm int}$ and for $\eps\in(0,1)$, let
$X^\eps$ be the biased voter model with generator (\ref{bvmgen}) and initial
state $x$. Define $(\mu^\eps_t)_{t\geq 0}$ and $(\mu_t)_{t\geq 0}$ as in
(\ref{mueps}) and (\ref{mu}). Then for each $0\leq t_1<\cdots<t_m$,
\be\label{fdd}
\P\big[(\mu^{\eps}_{t_1},\ldots,\mu^{\eps}_{t_m})\in\,\cdot\,\big]
\Ato{\eps}{0}\P\big[(\mu_{t_1},\ldots,\mu_{t_m})\in\,\cdot\,\big],
\ee
where $\Rightarrow$ denotes weak convergence of probability measures on
$\Mi(\R)^m$, and $\Mi(\R)$ is the space of locally finite measures on $\R$,
equipped with the topology of vague convergence.
\el

\bpro
Let
\be
\mu^{{\rm l},\eps}_t:=\!\!\!\!\sum_{i>L(X^\eps_{\eps^{-2}t})}\!\!\eps\de_{\eps i}
\quand
\mu^{{\rm r},\eps}_t:=\!\!\!\!\sum_{i>R(X^\eps_{\eps^{-2}t})}\!\!\eps\de_{\eps i}.
\ee
By Proposition~\ref{P:LRfdd} and the Skorohod
representation theorem \cite[Thm~6.7]{Bil99}, along any sequence $\eps_n\down
0$, we can couple our processes such that
\be
\eps_nL(X^{\eps_n}_{\eps_n^{-2}t_k})\asto{n}B_{t_k}
\quand
\eps_nR(X^{\eps_n}_{\eps_n^{-2}t_k})\asto{n}B_{t_k}
\quad{\rm a.s.}\qquad(1\leq k\leq m).
\ee
Then, for any continuous function $f:\R\to\R$ with compact support
\be\label{muvag}
\int_\R\!\mu^{{\rm l},\eps_n}_{t_k}(\di r)\,f(r)
\asto{n}\int_\R\!\mu_{t_k}(\di r)\,f(r)
\quad{\rm a.s.}\qquad(1\leq k\leq m),
\ee
and similarly for $\mu^{{\rm r},\eps_n}_{t_k}$, so using the fact that
$\mu^{{\rm r},\eps}_t\leq\mu^\eps_t\leq\mu^{{\rm l},\eps}_t$, we see that
(\ref{muvag}) holds with $\mu^{{\rm l},\eps_n}_{t_k}$ replaced by
$\mu^{\eps_n}_{t_k}$, first for $f\geq 0$ and then for general $f$ by
linearity. This proves that $\mu^{\eps_n}_{t_k}$ converges a.s.\ vaguely to
$\mu_{t_k}$ for each $1\leq k\leq m$. Since this holds for arbitrary $\eps_n\down
0$, (\ref{fdd}) follows.
\epro

\subsection{Tightness}\label{S:tight}

In this subsection, we complete the proof of Theorem~\ref{T:measure} by
showing tightness. 
Roughly speaking, we need to show that over a very short time interval, 
there cannot be too many 0's invading far into the region dominated by 1's 
and vice versa. 
The first scenario will be ruled out by a direct comparison with the unbiased 
voter model, for which tightness was proved in \cite{AS11}. 
To rule out the second scenario, the key observation is that once enough 1's 
invade far into the region dominated by 0's, they will persist, leading to a large 
displacement of the interface at the end of the short time interval. 

We let $\li\mu,\phi\re:=\int \phi\,\di\mu$ denote the
integral of a function $\phi$ with respect to a measure $\mu$, and we write
$\Ci^2_{\rm c}(\R)$ for the space of compactly supported, twice continuously
differentiable functions $f:\R\to\R$.
Let $\Ki$ denote the space of all measures $\mu\in\Mi(\R)$ such that
$\mu([-n,n])\leq 2n+1$ for $n=1,2,\ldots$. Then $\Ki$ is a compact subset of
$\Mi(\R)$ and $\mu^\eps_t\in\Ki$ for all $\eps\in[0,1)$ and $t\geq 0$. In view
of this, by Jakubowski's tightness criterion \cite[Thm 3.6.4]{DA93}, for given
$\eps_n\down 0$, the laws of $\{(\mu^{\eps_n}_t)_{t\geq 0}\}_{n\in\N}$ are tight on
$\Di([0,\infty),\Mi(\R))$ if and only if
\begin{itemize}
\item[{\bf(J)}] (Tightness of evaluations) For each $f\in\Ci^2_{\rm c}(\R)$,
the laws of $\{(\li\mu^{\eps_n}_t,f\re)_{t\geq 0}\}_{n\in\N}$ are tight on
$\Di([0,\infty),\R)$.
\end{itemize}
To verify tightness of the laws of the real-valued processes
$(\li\mu^{\eps_n}_t,f\re)_{t\geq 0}$, we will use the following lemma.

\bl[Tightness criterion]\label{L:ticri}
Let\label{L:tight} $\xi^n\in\Di(\half,\R)$ and assume that
\be\label{tight1}
\lim_{\de\down 0}\sum_{i=0}^{\lfloor T/\de\rfloor}
 \limsup_{n\to\infty}\P^x\big[\sup_{t\in[i\de,(i+1)\de)}
 |\xi^n_t-\xi^n_{i\de}|\geq\eta\big]=0\qquad(\eta>0,\ T<\infty).
\ee
Then the laws $\P[\xi^n\in\,\cdot\,]$ are tight on $\Di(\half,\R)$ and each
weak limit point is concentrated on $\Ci(\half,\R)$.
\el

\bpro
It is well known \cite[Thm~15.5]{Bil99} that the conclusion of the lemma is implied by
\be\label{tight3}
\lim_{\de\down 0}\limsup_{n\to\infty}
\P^x\big[\sup_{0\leq s<t\leq T:\,|s-t|\leq\de} 
|\xi^n_t-\xi^n_s|\geq\eta\big]=0\qquad(\eta>0,\ T<\infty).
\ee
If $|\xi^n_t-\xi^n_s|\geq\eta$ for some $0\leq s<t\leq T$ with $|s-t|\leq\de$,
then there must exist $0\leq i\leq\lfloor T/\de\rfloor$ and $i\de\leq
s<t<(i+1)\de$ such that $|\xi^n_t-\xi^n_s|\geq\eta/2$, and hence
$\sup_{t\in[i\de,(i+1)\de)}|\xi^n_t-\xi^n_{i\de}|\geq\eta/4$. This shows that
(\ref{tight1}) implies (\ref{tight3}).
\epro

%We assume, for the moment, that the kernel $a$ is non-nearest-neighbor.
We will establish tightness for $\{(\mu^{\eps_n}_t)_{t\geq0}\}_{n\in\N}$ by
delicate comparisons between biased and unbiased voter models using results
from \cite{AS11} that we now cite.

%Although it is not formulated as a lemma there, we essentially cite the
%following result from \cite{AS11}.

\bl[Continuity estimate for the unbiased model]
Let\label{L:unbias} $\P^x$ denote the law of the unbiased voter model
$(X^0_t)_{t\geq 0}$ started in $X^0_0=x$ and let
\be\label{nudef}
\nu^\eps_t:=\sum_{i\in\Z}\eps X^0_{\eps^{-2}t}(i)\de_{\eps i}
\qquad(\eps>0,\ t\geq 0).
\ee
Then for each $f\in\Ci^2_{\rm c}(\R)$, there exist $C<\infty$ and $t_0,\eps_0>0$
such that for all $0\leq t\leq t_0$ and $0<\eps\leq\eps_0$,
\be\label{unbias}
\P^x\big[\big|\li\nu^\eps_t,f\re-\li\nu^\eps_0,f\re\big|\geq\de\big]
\leq C\de^{-2}t^{1/4}\qquad(x\in\{0,1\}^\Z,\ \de>0).
\ee
\el

\bpro
This is proved in Section~2.1 of \cite{AS11}, as a first step towards proving
that the laws of the processes in (\ref{nudef}) are tight. The proof uses the
duality between the voter model and coalescing random walks to derive
estimates for the mean and variance of $\li\nu^\eps_t,f\re$. Crucially, the
bounds are independent of the initial state $x$ and only assume the
properties (i)--(iii) of the kernel $a(\,\cdot\,)$ that we also use.
\epro

Recall that the main result of \cite{AS11} is that if an unbiased voter model
$(X^0_t)_{t\geq0}$ is started in the Heaviside state $x_{\rm hv}$, then
$(\nu^\eps_t)_{t\geq 0}$ converges to $(1_{\{y\geq\ti B_t\}}\di y)_{t\geq 0}$ as
$\eps\down0$, where $\ti B_t:=W_{\sig^2 t}$ is a Brownian motion with
diffusion coefficient $\sig^2$. Indeed, their proof can be extended to
more general initial configurations.

%that are random, but converge to $1_{\{y\geq 0\}}\di y$.

\bl[Invariance principle for the unbiased voter model]
Let\label{L:vmeasure} $\eps_n\down0$. For each $n$, let $(X^{0,n}_t)_{t\geq 0}$
be an unbiased voter model started in a deterministic initial state and define
$\nu^{\eps_n}_t$ as in (\ref{nudef}) but with $X^0$ replaced by
$X^{0,n}$. Assume that $\nu^{\eps_n}_0$ converges in the vague topology to
$1_{\{y\geq 0\}}\di y$ as $n\to\infty$. Then
\be
\P\big[(\nu^{\eps_n}_t)_{t\geq0}\in\,\cdot\,\big]
\Ato{n}{\infty}\P\big[(\nu_t)_{t\geq0}\in\,\cdot\,\big],
\ee
where $\nu_t:=1_{\{y\geq \ti B_t \}}\di y$ is the Lebesgue measure on a half
line whose boundary is given by the Brownian motion $(\ti B_t)_{t\geq 0}$ with
diffusion coefficient $\sig^2$.
\el

\bpro
When the initial configuration is $x_{\rm hv}$, tightness and convergence in
finite dimensional distributions are shown in Sections~2.1 and 2.2 of
\cite{AS11}, respectively. To prove tightness, \cite{AS11} also use
Jakubowski's tightness criterion, where Aldous' tightness criterion (see
e.g.\ \cite[Thm 1]{Ald78}) is applied to verify {\bf(J)} that is given above Lemma~\ref{L:ticri}. 
A crucial ingredient is the estimate (\ref{unbias}), which holds for general initial
configurations. In view of this, their proof of tightness holds regardless of
the initial condition.

The proof of convergence of the finite dimensional distributions is based on a
first and second moment calculation using duality and the fact that a
collection of dual coalescing random walks converges to a collection of
coalescing Brownian motions. For this part of the argument, it suffices if the
initial configuration $\nu^{\eps_n}_0$ converges in the vague topology to
$1_{\{y\geq 0\}}\di y$ as $n\to\infty$.
\epro

We now prove tightness by verifying Jakubowski's tightness criterion {\bf(J)},
using Lemmas~\ref{L:tight}, \ref{L:unbias}, and \ref{L:vmeasure}, and
judicious comparisons between biased and unbiased voter models, thereby
completing the proof of our main result Theorem~\ref{T:measure}.%\med

\bpro[Proof~of~Theorem~\ref{T:measure}]
If $a$ is the nearest-neighbor kernel $a(-1)=\ha=a(1)$, then the Heaviside
state is a trap for the process modulo translations.  As in (\ref{tau0}), let
$\tau^\eps_0$ denote the trapping time. It has been shown below (\ref{asfin})
that in this case, the biased voter model observed until $\tau^\eps_0$
converges in law to the unbiased voter model observed until $\tau^0_0$, and
that $\tau^0_0$ is finite a.s. In view of this, in this case,
Theorem~\ref{T:measure} follows trivially from Theorem~\ref{T:midconv}. We
assume therefore without loss of generality that $a$ is non-nearest-neighbor.

Convergence of finite dimensional distributions has already been proved
in Lemma~\ref{L:fdd}, so it suffices to show tightness. As argued at the
beginning of this section, by Jakubowski's tightness criterion, it suffices to
show that for each $f\in\Ci^2_{\rm c}(\R)$, the laws of the processes
$\{(\li\mu^{\eps_n}_t,f\re)_{t\geq 0}\}_{n\in\N}$ are tight on
$\Di([0,\infty),\R)$ along any sequence $\eps_n\down 0$.
By linearity, it suffices to consider nonnegative $f$.
We fix $f\geq 0$ and apply Lemma~\ref{L:tight} to the real-valued processes
$(\li\mu^{\eps_n}_t,f\re)_{t\geq 0}$. We fix $\eta>0$ and for each $n$ and
$s\geq 0$ define
\bc
\dis\tau^{n,+}_s&:=&\dis\inf\big\{t\geq 0:
\li\mu^{\eps_n}_{s+t},f\re-\li\mu^{\eps_n}_s,f\re\geq\eta\big\},\\[5pt]
\dis\tau^{n,-}_s&:=&\dis\inf\big\{t\geq 0:
\li\mu^{\eps_n}_{s+t},f\re-\li\mu^{\eps_n}_s,f\re\leq-\eta\big\}.
\ec
We will prove the lower and upper bounds
\be\left.\ba{r@{\ }l}\label{tight}
{\rm(i)}&\dis\lim_{\de\down 0}\de^{-1}\sup_{s\geq 0}
 \limsup_{n\to\infty}\P[\tau^{n,+}_s<\de]=0,\\[5pt]
{\rm(ii)}&\dis\lim_{\de\down 0}\de^{-1}\sup_{s\geq 0}
 \limsup_{n\to\infty}\P[\tau^{n,-}_s<\de]=0
\ea\right\}\qquad(\eta>0),
\ee
which together imply (\ref{tight1}) and hence tightness for the laws of
$(\li\mu^{\eps_n}_t,f\re)_{t\geq 0}$.

To prove (\ref{tight}) (i), we note that if $\tau^{n,+}_s<\de$ and the increment of the process
$(\li\mu^{\eps_n}_t,f\re)_{t\geq 0}$ during the remaining time from $s+\tau^{n,+}_s$ to $s+\delta$ 
is no less than $-\eta/2$,  then the increment of the process between $s$ and $s+\delta$
is larger than $\eta/2$.
Thus,
\be\label{Cden}
\P\big[\li\mu^{\eps_n}_{s+\de},f\re-\li\mu^{\eps_n}_s,f\re\geq\eta/2\big]
\geq C_{\de,n}\P[\tau^{n,+}_s<\de]
\ee
where
\be\label{Cdef}
C_{\de,n}:=\inf_{0\leq t<\de}\inf_x
\P^x\big[\li\mu^{\eps_n}_t,f\re-\li\mu^{\eps_n}_0,f\re\geq-\eta/2\big],
\ee
and we have conditioned on $\tau=\tau^{n,+}_s$ and $x=X^{\eps_n}_{s+\tau}$
and used the strong Markov property. We couple the biased voter model
started in $X^{\eps_n}_0=x$ to an unbiased voter model started in $X^0_0=x$
in such a way that $X^{\eps_n}_t\geq X^0_t$ for all $t\geq 0$. Defining
$\nu^{\eps_n}_t$ as in (\ref{nudef}), using that $f\geq 0$, it follows that
\be\label{downcoup}
\P^x\big[\li\mu^{\eps_n}_t,f\re-\li\mu^{\eps_n}_0,f\re\geq-\eta/2\big]
\geq\P^x\big[\li\nu^{\eps_n}_t,f\re-\li\nu^{\eps_n}_0,f\re\geq-\eta/2\big].
\ee
Using the symmetry between zeros and ones in the unbiased voter model,
we obtain by Lemma~\ref{L:unbias} that $C_{\de,n}\geq
1-C\eta^{-2}\de^{1/4}\geq 1/2$ for all $\de$ small enough and $n$ large
enough. Inserting this into (\ref{Cden}) and using the convergence of the
finite dimensional distributions (Lemma~\ref{L:fdd}), we find that for $\de$
small enough,
\bc
\dis\limsup_{n\to\infty}\P[\tau^{n,+}_s<\de]
&\leq&\dis 2\limsup_{n\to\infty}
\P\big[\li\mu^{\eps_n}_{s+\de},f\re-\li\mu^{\eps_n}_s,f\re\geq\eta/2\big]\\[5pt]
&=&\dis 2
\P\big[\int f(x)1_{\{x\geq B_{s+\de}\}}\di x
-\int f(x)1_{\{x\geq B_s\}}\di x\geq\eta/2\big]\\[5pt]
&\leq&\dis 2
\P\big[|B_{s+\de}-B_s|\geq\ha\eta\|f\|_\infty\big],
\ec
where $B_t=W_{\sig^2t}-\ha\sig^2t$ is a Brownian motion with drift
$-\ha\sig^2$ and diffusion constant $\sig^2$. It is easy to see the right-hand
side is $o(\de)$, uniformly in $s$, proving (\ref{tight}) (i).

The argument for (\ref{tight}) (ii) is similar, but not quite the same.
In this case, we couple biased and unbiased voter models started in the same
initial state at time $s$ to bound
\be\label{taus}
\P[\tau^{n,-}_s<\de]\leq\P[\sig^{n,-}_s<\de],
\ee
where
\be
\sig^{n,-}_s:=\inf\big\{t\geq 0:
\li\nu^{\eps_n}_{s+t},f\re-\li\nu^{\eps_n}_s,f\re\leq-\eta\big\}.
\ee
and we use that $f\geq 0$. Arguing as in (\ref{Cden}) and (\ref{Cdef}), applying
Lemma~\ref{L:unbias} directly without the need of the coupling in
(\ref{downcoup}), allows us to estimate, for $\de$ small enough,
\bc
\dis\limsup_{n\to\infty}\P[\sig^{n,-}_s<\de]
&\leq&\dis 2\limsup_{n\to\infty}
\P\big[\li\nu^{\eps_n}_{s+\de},f\re-\li\nu^{\eps_n}_s,f\re\geq\eta/2\big]\\[5pt]
&=&\dis 2
\P\big[\int f(x)1_{\{x\geq\ti B_{s+\de}\}}\di x
-\int f(x)1_{\{x\geq\ti B_s\}}\di x\geq\eta/2\big]\\[5pt]
&\leq&\dis 2
\P\big[|\ti B_{s+\de}-\ti B_s|\geq\ha\eta\|f\|_\infty\big],
\ec
where we have used Lemma~\ref{L:vmeasure} instead of Lemma~\ref{L:fdd} and
$(\ti B_t)_{t\geq s}$ is a Brownian motion with zero drift and diffusion
constant $\sig^2$, started at time $s$ in $\ti B_s=B_s$. Again, the right-hand
side is $o(\de)$, which together with (\ref{taus}) proves (\ref{tight}) (ii).
\epro

%hiero

\appendix

%The flow paper (AMS) contains a useful appendix on the Hausdorff topology.
%Etheridge, Freeman and Straulino [34] work with compact sets of cadlag paths
%My lecture notes with Anita, Thm 1.8 is a usueful ref for Skorohod topology.
%Otherwise, Jacod & shiryayev, Ethier & Kurtz are basic references.

\section{Appendix}

\subsection{Locally uniform convergence}

For any metrizable space $E$, we let $\Di(\half,E)$ denote the space of
\cadlag\ functions (i.e., right-continuous functions with left limits)
$w:\half\to E$, equipped with the Skorohod topology \cite{EK86,Bil99}, and we
let $\Ci(\half,E)$ denote the subspace of continuous functions. It is
well-known that $\Di(\half,E)$ is Polish if $E$ is
\cite[Thm~3.5.6]{EK86}. Moreover, a sequence $w_n\in\Di(\half,E)$ converges to
a limit $w\in\Ci(\half,E)$ if and only if $w_n\to w$ locally uniformly on
compact sets \cite[Lemma~3.10.1]{EK86}.  We recall the following well-known lemma.

\bl[Convergence criterion]
Let\label{L:cc} $E$ be a metrizable space and let $d$ be any metric generating
the topology on $E$. Let $w_n,w:\half\to E$ be functions and assume that $w$
is continuous. Then $w_n\to w$ locally uniformly if and only if $w_n(t_n)\to
w(t)$ for all $t_n,t\geq 0$ such that $t_n\to t$.
\el

It is not hard to see that locally uniform convergence of functions implies
locally uniform convergence of their compositions and inverses. Moreover, for
monotone functions, pointwise convergence is equivalent to locally uniform
convergence. %For completeness, we provide proofs.

\bl[Convergence of composed functions]
Let\label{L:compconv} $E$ be a metrizable space, let
$\la_n,\la:\half\to\half$ be nondecreasing, and let
$w_n,w:\half\to E$. Assume moreover that $\la,w$ are continuous.
Then $\la_n\to\la$ locally uniformly and $w_n\to w$ locally uniformly imply that
$w_n\circ\la_n\to w\circ\la$ locally uniformly.
\el

\bpro
By Lemma~\ref{L:cc} $t_n\to t$ implies $\la_n(t_n)\to\la(t)$ and hence
$w_n(\la(t_n))\to w(\la(t))$. Since this holds for general $t_n\to t$,
the claim now follows from Lemma~\ref{L:cc}.
\epro

\bl[Convergence of nondecreasing functions]
Let\label{L:monpt} $\la_n,\la:\half\to\half$ be nondecreasing, and assume
that $\la$ is continuous. Let $D\sub\half$ be dense. Then $\la_n\to\la$
locally uniformly if and only if for all $t\in D$ there exist $t_n\geq 0$
such that $t_n\to t$ and $\la_n(t_n)\to\la(t)$.
\el

\bpro
The necessity of the condition is clear. To prove the sufficiency, by
Lemma~\ref{L:cc} it suffices to show that $t_n\to t$ implies
$\la_n(t_n)\to\la(t)$. Fix $t^\pm\in D$ with $t^-<t<t^+$ and choose
$t^\pm_n\geq 0$ such that $t^\pm_n\to t^\pm$ and
$\la_n(t^\pm_n)\to\la(t^\pm)$. Then $t^-_n<t_n<t^+_n$ for $n$ sufficiently
large, and hence, since the $\la_n$ are nondecreasing,
$\la_n(t^-_n)\leq\la_n(t_n)\leq\la_n(t^+_n)$ for $n$ sufficiently large. It
follows that $\la(t_-)\leq\liminf_{n\to\infty}\la_n(t_n)$ and
$\limsup_{n\to\infty}\la_n(t_n)\leq\la(t_+)$.
Using the density of $D$ and the continuity of $\la$, we conclude that
$\la_n(t_n)\to\la(t)$.
\epro

For any $\la\in\Ci(\half,\half)$, let $\la(\half):=\{\la(t):t\in\half\}$
denote the image of $\half$ under $\la$. If $\la$ is strictly increasing and
$\la(\half)=\half$, then $\la$ has an inverse $\la^{-1}$.

\bl[Convergence of inverse functions]
Let\label{L:invconv} $\la_n,\la\in\Ci(\half,\half)$ be strictly increasing with
$\la_n(\half)=\half$ and $\la(\half)=\half$. Then $\la_n\to\la$
locally uniformly if and only if $\la^{-1}_n\to\la^{-1}$ locally uniformly.
\el

\bpro
Let $G:=\{(t,\la(t)):t\geq 0\}$ denote the graph of $\la$
and similarly, let $G_n$ denote the graph of $\la_n$. Then
the graph of $\la^{-1}$ is $G^{-1}=\{(\la(t),t):t\geq 0\}$ and similarly for
the graph $G_n^{-1}$ of $\la_n^{-1}$. Lemma~\ref{L:monpt} tells us that
$\la_n\to\la$ locally uniformly if and only if for all $(t,s)\in G$ there
exist $(t_n,s_n)\in G_n$ such that $(t_n,s_n)\to(t,s)$. Clearly, this holds if
and only if $G_n^{-1}$ and $G^{-1}$ satisfy the same condition, which is equivalent
to $\la^{-1}_n\to\la^{-1}$ locally uniformly.
\epro

Let\label{L:top} $X_n$ be random variables taking values in a Polish space
$E$, and let $x\in E$. Then it is not hard to see that the following
statements are equivalent
\begin{enumerate}
\item $\dis\P[X_n\in\,\cdot\,]\Asto{n}\de_x$,
\item $\dis\P\big[X_n\not\in A\big]\asto{n}0$ for all $A\in\Ni_x$,
\end{enumerate}
where $\Rightarrow$ denotes weak convergence of probability measures and
$\Ni_x$ is a fundamental system of neighborhoods of $x$. If these conditions
are fulfilled, then we say that the $X_n$ converge to $x$ \emph{in
  probability} and denote this as $X_n\toP x$. In particular, let $E$ be a
Polish space and $d$ a metric generating the topology on $E$, let $W_n$ be
random variables with values in $\Di_E\half$, and let $w\in\Ci_E\half$. Then
$W_n\toP w$ with respect to the Skorohod topology if and only if
\be
\sup_{0\leq t\leq T}d\big(W_n(t),w(t)\big)\atoP{n}{\infty}0\qquad(T<\infty).
\ee

Because we will need these in our proofs, for completeness, we provide proofs
for two additional simple lemmas which lift Lemmas~\ref{L:compconv} and
\ref{L:monpt} to convergence in law and in probability, respectively.

\bl[Convergence of time-changed processes]
Let\label{L:timconv} $Y=(Y_t)_{t\geq 0}$ and $Y^n=(Y^n_t)_{t\geq 0}$ be
stochastic processes with \cadlag\ sample paths, taking values in a Polish
space $E$. Let $S=(S_t)_{t\geq 0}$ and $S^n=(S^n_t)_{t\geq 0}$ be real-valued
stochastic processes whose sample paths are \cadlag, nondecreasing, and
satisfy $S_0=0$ resp.\ $S^n_0=0$. Assume moreover that $Y=(Y_t)_{t\geq 0}$
and $S=(S_t)_{t\geq 0}$ have continuous sample paths, and that
\be
\P[(Y^n_t,S^n_t)_{t\geq 0}\in\,\cdot\,]
\Asto{n}\P[(Y_t,S_t)_{t\geq 0}\in\,\cdot\,],
\ee
where $\Rightarrow$ denotes weak convergence with respect to the Skorohod
topology. Then
\be\label{timconv}
\P\big[(Y^n_{S^n_t})_{t\geq 0}\in\,\cdot\,\big]
\Asto{n}
\P\big[(Y_{S_t})_{t\geq 0}\in\,\cdot\,\big].
\ee
\el

\bpro
By the Skorohod representation theorem \cite[Thm~6.7]{Bil99}, we can couple
our random variables such that $(Y^n_t,S^n_t)_{t\geq 0}$ converges a.s.\ to
$(Y_t,S_t)_{t\geq 0}$ with respect to the Skorohod topology. Now
Lemma~\ref{L:compconv} implies that $(Y^n_{S^n_t})_{t\geq 0}$ converges a.s.\ to
$(Y_{S_t})_{t\geq 0}$ w.r.t.\ the same topology, and hence (\ref{timconv})
follows.
\epro

\bl[Convergence of nondecreasing functions]
Let\label{L:pmonpt} $S^n=(S^n_t)_{t\geq 0}$ be real-valued
stochastic processes whose sample paths are \cadlag, nondecreasing, and
satisfy $S_0=0$ resp.\ $S^n_0=0$. Let $\la:\half\to\half$ be continuous.
Then the following statements are equivalent:
\[
{\rm(i)}\quad\sup_{0\leq t\leq T}\big|S^n_t-\la_t\big|\atoP{n}{\infty}0
        \qquad(T<\infty),\qquad
{\rm(ii)}\quad S^n_t\atoP{n}{\infty}\la_t\qquad(t\geq 0).
\]
\el

\bpro
The implication (i)$\volgt$(ii) is trivial. To prove the converse, let
$\{t_k:k\in\N\}$ be countable and dense. Then
\be
\sup_{0\leq k\leq m}|S^n_{t_k}-\la_{t_k}|\atoP{n}{\infty}0\qquad(m<\infty),
\ee
which says that the process $k\mapsto S^n_{t_k}$ converges in probability
to $k\mapsto\la_{t_k}$ with respect to the product topology on $\R^\N$.
By the Skorohod representation theorem \cite[Thm~6.7]{Bil99}, we can couple
our random variables such that
\be
S^n_{t_k}\asto{n}\la_{t_k}\quad{\rm a.s.}\qquad(k\in\N).
\ee
By Lemma~\ref{L:monpt}, it follows that $\sup_{0\leq t\leq
  T}\big|S^n_t-\la_t\big|$ converges a.s.\ to zero for all $T<\infty$,
which implies (i).
\epro

\subsection{A weak law of large numbers}

In this subsection we prove two simple versions of the weak law of large
numbers. Lemma~\ref{L:FLLN} below is used in the proof of
Theorem~\ref{T:midconv}. The following lemma would be completely
standard if the law of $(V_{n,i})_{i\geq 1}$ would not depend on $n$.

\bl[A weak law of large numbers]
For\label{L:LLN} each $n\geq 1$, let $(V_{n,i})_{i\geq 1}$ be
i.i.d.\ nonnegative random variables, and let $m_n\geq 1$ be integers
such that $\lim_{n\to\infty}m_n=\infty$. Assume that
\be\label{trunc}
\sup_{n\geq 1}\E\big[|V_{n,1}|\big]<\infty\quand
\E\big[|V_{n,1}|;|V_{n,1}|>t m_n]\asto{n}0\qquad(t>0).
\ee
Then
\be
\E[V_{n,1}]-\frac{1}{m_n}\sum_{i=1}^{m_n}V_{n,i}\atoP{n}{\infty}0,
\ee
where $\overset{\rm P}{\rightarrow}$ denotes convergence in probability.
\el

\bpro
Define truncated random variables by $\ov V_{n,i}:=V_{n,i}1_{\{|V_{n,i}|\leq
  m_n\}}$. Then (\ref{trunc}) implies that
\be
\P\Big[\sum_{i=1}^{m_n}\ov V_{n,i}\neq\sum_{i=1}^{m_n}V_{n,i}\Big]
\leq m_n\P\big[|V_{n,1}|>m_n\big]\leq\E\big[|V_{n,1}|;|V_{n,1}|>m_n]\asto{n}0.
\ee
Since (\ref{trunc}) moreover implies that
\be\label{wlln3}
\big|\E[\ov V_{n,1}]-\E[V_{n,1}]\big|
\leq\E[|V_{n,1}|;|V_{n,1}|>m_n]\ato{n}{\infty}0,
\ee
it suffices to prove the statement with $V_{n,i}$ replaced by $\ov V_{n,i}$.
For any $\de>0$, Chebyshev's inequality gives
\be\label{Cheb}
P\Big[\Big|\frac{1}{m_n}\sum_{i=1}^{m_n}\ov V_{n,i}-\E[\ov V_{n,1}]\Big|>\de\Big]
\leq\de^{-2}\frac{1}{m_n}\var(\ov V_{n,1}).
\ee
Using the layer cake representation, we estimate
\be
\var(\ov V_{n,1})\leq\E[\ov V^2_{n,1}]
\leq2\int_0^{m_n}\!x\,\P\big[|V_{n,1}|>x\big]\di x
\leq2\int_0^{m_n}\E\big[|V_{n,1}|;|V_{n,1}|>x\big]\di x.
\ee
It follows that the right-hand side of (\ref{Cheb}) can be estimated
by
\be
2\de^{-2}\int_0^1\E\big[|V_{n,1}|;|V_{n,1}|>tm_n\big]\di t,
\ee
which tends to zero by (\ref{trunc}), using dominated convergence.
\epro

%\bl[Functional law of large numbers]
%For\label{L:FLLN} each $n\geq 1$, let $(V_{n,i})_{i\geq 1}$ be
%i.i.d.\ nonnegative random variables, and let $\eps_n>0$ be constants
%such that $\lim_{n\to\infty}\eps_n=0$. Assume that
%$\sup_n\E[|V_{n,1}|]<\infty$,
%\be\label{trunc2}
%\lim_{n\to\infty}\E[V_{n,1}]=c<\infty,
%\quand\E\big[|V_{n,1}|;|V_{n,1}|>t/\eps_n]\asto{n}0\qquad(t>0).
%\ee
%Define $f_n:\half\to\half$ by
%\be
%f_n(t):=\eps_n\sum_{i=1}^{\lfloor\eps_n^{-1}t\rfloor}V_{n,i}\qquad(t\geq 0).
%\ee
%Then
%\be
%\sup_{0\leq t\leq T}\big|ct-f_n(t)\big|\atoP{n}{\infty}0\qquad(T>0),
%\ee
%where $\overset{\rm P}{\rightarrow}$ denotes convergence in probability.
%\el
%
%\bpro
%Applying Lemma~\ref{L:LLN} to $m_n=\lfloor\eps_n^{-1}t\rfloor$, we see
%that $f_n(t)$ converges in probability to $ct$ for each fixed $t>0$.
%If the $V_{n,i}$ are nonnegative, the claim now follows from
%Lemma~\ref{L:pmonpt}.
%\epro

\bl[Functional law of large numbers]
For\label{L:FLLN} each $n\geq 1$, let $(V_{n,i})_{i\geq 1}$ be
i.i.d.\ nonnegative random variables, and let $\eps_n>0$ be constants
such that $\lim_{n\to\infty}\eps_n=0$. Assume that
\be\label{trunc2}
\lim_{n\to\infty}\E[V_{n,1}]=c<\infty
\quand\E\big[V_{n,1};V_{n,1}>t/\eps_n]\asto{n}0\qquad(t>0).
\ee
Define $f_n:\half\to\half$ by
\be
f_n(t):=\eps_n\sum_{i=1}^{\lfloor\eps_n^{-1}t\rfloor}V_{n,i}\qquad(t\geq 0).
\ee
Then
\be
\sup_{0\leq t\leq T}\big|ct-f_n(t)\big|\atoP{n}{\infty}0\qquad(T>0),
\ee
where $\overset{\rm P}{\rightarrow}$ denotes convergence in probability.
\el

\bpro
Since the $V_{n,i}$ are nonnegative, the condition
$\lim_{n\to\infty}\E[V_{n,1}]=c<\infty$ implies that $\sup_{n\geq
n_0}\E\big[V_{n,1}\big]<\infty$ for $n_0$ sufficiently large. Applying
Lemma~\ref{L:LLN} to $m_n=\lfloor\eps_n^{-1}t\rfloor$, we see that $f_n(t)$
converges in probability to $ct$ for each fixed $t>0$. The claim now follows
from Lemma~\ref{L:pmonpt}.
\epro

\subsection{Uniform ergodicity}

The following lemma, which we apply in Subsection~\ref{S:fdd}, gives
sufficient conditions for the speed of convergence to equilibrium to be
uniform for a sequence of continuous-time Markov chains.

\bl[Uniform ergodicity]
Let\label{L:unerg} $S$ be a countable set and for each $n\in\N\cup\{\infty\}$,
let $X^n=(X^n_t)_{t\geq 0}$ be a positive recurrent, irreducible continuous-time
Markov chain with state space $S$ and invariant law $\pi_n$. Assume that as
$n\to\infty$, the jump rates of $X^n$ converge pointwise to the jump rates of
$X^\infty$, and the invariant laws $\pi_n$ converge weakly to
$\pi_\infty$. Then, for each $x\in S$,
\be\label{unerg}
\sup_{n\in\N\cup\{\infty\}}\big\|\P^x[X^n_t\in\,\cdot\,]-\pi_n\big\|
\asto{t}0,
\ee
where $\|\,\cdot\,\|$ denotes the total variation norm.
\el

\bpro
Fix $z\in S$. Let $(X^n_t)_{t\geq 0}$ and $(\ti X^n_t)_{t\geq 0}$ be
independent processes with the same jump rates and let
$\tau^n_{(z,z)}:=\inf\{t\geq 0:X_t=z=\ti X_t\}$. 
Using Doeblin's coupling, we can couple two
processes by declaring them to be equal after $\tau^n_{(z,z)}$.
Therefore, we see that
\be
\big\|\P^x[X^n_t\in\,\cdot\,]-\P^y[X^n_t\in\,\cdot\,]\big\|
\leq\P^{(x,y)}[t<\tau^n_{(z,z)}]\qquad(x,y\in S,\ t\geq 0),
\ee
and hence
\bc\label{ncoup}
\dis\big\|\P^x[X^n_t\in\,\cdot\,]-\pi_n\big\|
&=&\dis\big\|\sum_{y\in S}\pi_n(y)\big(\P^x[X^n_t\in\,\cdot\,]
-\P^y[X^n_t\in\,\cdot\,]\big)\big\|\\[5pt]
&\leq&\dis\sum_{y\in S}\pi_n(y)\P^{(x,y)}[t<\tau^n_{(z,z)}]
\qquad(x\in S,\ t\geq 0).
\ec
Since the jump rates converge, the probability $\P^{(x,y)}[t<\tau^n_{(z,z)}]$
converges pointwise as $n\to\infty$ for each $y\in S$.
Using also that $\pi_n\Rightarrow\pi_\infty$, which implies that the measures
$\pi_n$ are tight, this is easily seen to imply that the right-hand side of
(\ref{ncoup}) converges and hence
\be\label{alm}
\limsup_{n\to\infty}\big\|\P^x[X^n_t\in\,\cdot\,]-\pi_n\big\|
\leq\sum_{y\in S}\pi_\infty(y)\P^{(x,y)}[t<\tau^\infty_{(z,z)}]
\qquad(x\in S,\ t\geq 0).
%\asto{t}0\qquad(x\in S),
\ee
The joint process $(X^\infty_t,\ti X^\infty_t)_{t\geq 0}$ is irreducible and
has an invariant law $\pi_\infty\otimes\pi_\infty$, which implies positive
recurrence. In view of this, the right-hand side of (\ref{alm}) converges to
zero as $t\to\infty$ for each fixed $x\in S$. Since $X^n$ is positive
recurrent and hence ergodic for each $n\in\N$ and since the total variation
distance to the invariant measure is a nonincreasing function of time,
\bc
\dis\limsup_{t\to\infty}\,
\sup_n\big\|\P^x[X^n_t\in\,\cdot\,]-\pi_n\big\|
&\leq&\dis\limsup_{t\to\infty}\,\sup_{n\geq N}
\big\|\P^x[X^n_t\in\,\cdot\,]-\pi_n\big\| \\ [5pt]
&\leq&\dis\sup_{n\geq N}\big\|\P^x[X^n_T\in\,\cdot\,]-\pi_n\big\|
%\qquad(T<\infty),
\ec
for each $N,T<\infty$,
%\be\ba{l}
%\dis\limsup_{t\to\infty}
%\sup_{n\in\N\cup\{\infty\}}\big\|\P^x[X^n_t\in\,\cdot\,]-\pi_n\big\|\\[5pt]
%\dis\quad
%\leq\limsup_{t\to\infty}\sup_{n<N}\big\|\P^x[X^n_t\in\,\cdot\,]-\pi_n\big\|
%\vee\limsup_{t\to\infty}\sup_{n\geq N}
%\big\|\P^x[X^n_t\in\,\cdot\,]-\pi_n\big\|\\[5pt]
%\dis\quad\leq\sup_{n\geq N}\big\|\P^x[X^n_T\in\,\cdot\,]-\pi_n\big\|
%\qquad(T<\infty),
%\ec
where in view of (\ref{alm}) the right-hand side can be made arbitrary small
by choosing $N$ and $T$ large enough.
\epro

%Indeed, it follows that for each $\eps>0$, there exists an $n_0$ such that
%\be
%\sup_{n\geq n_0}\big\|\P^x[X^n_t\in\,\cdot\,]-\pi_n\big\|
%\leq\eps+\sum_{y\in S}\pi_\infty(y)\P^{(x,y)}[t<\tau^\infty_{(z,z)}],
%\ee
%where the supremum runs over all $n\in\N\cup\{\infty\}$ such that $n\geq n_0$.
%Then
%\be
%\sup_{n\in\N\cup\{\infty\}}\big\|\P^x[X^n_t\in\,\cdot\,]-\pi_n\big\|
%=\sup_{n<n_0}\big\|\P^x[X^n_t\in\,\cdot\,]-\pi_n\big\|
%\vee\sup_{n\geq n_0}\big\|\P^x[X^n_t\in\,\cdot\,]-\pi_n\big\|,
%\ee
%so, using the fact that $X^n$ is positive recurrent and hence ergodic for each
%$n$, we obtain
%\be
%\limsup_{t\to\infty}\sup_{n\in\N\cup\{\infty\}}
%\big\|\P^x[X^n_t\in\,\cdot\,]-\pi_n\big\|\leq\eps
%\ee
%for all $\eps>0$.

\end{document}